\documentclass[11 pt]{amsart}
\usepackage{amssymb}
\usepackage{amscd}

\newtheorem{theorem}{Theorem}
\newtheorem{corollary}[theorem]{Corollary}
\newtheorem{lemma}[theorem]{Lemma}
\newtheorem{proposition}[theorem]{Proposition}

\theoremstyle{definition}
\newtheorem{definition}{Definition}

\numberwithin{equation}{section}

\newcommand{\ba}{\mathbb{A}}
\newcommand{\bc}{\mathbb{C}}
\newcommand{\bp}{\mathbb{ P}}
\newcommand{\bz}{\mathbb{Z}}

\newcommand{\bq}{\mathbb{Q}}

\newcommand{\cc}{\mathcal{C}}

\newcommand{\cl}{\mathcal{L}}
\newcommand{\cs}{\mathcal{S}}

\newcommand{\ch}{\mathcal{H}}
\newcommand{\cm}{\mathcal{M}}

\newcommand{\ca}{\mathcal{A}}

\newcommand{\bg}{\bigskip}
\newcommand{\med}{\medskip}
\newcommand{\la}{\longrightarrow}
\newcommand{\bfl}{\begin{flushleft}}
\newcommand{\efl}{\end{flushleft}}

\newcommand{\om}{\Omega}
 \newcommand{\kh}{K_{hol}}
\newcommand{\kt}{K_{top}}
\newcommand{\ka}{K_{alg}}
\newcommand{\tkh}{\tilde \kh}

\newcommand{\tkt}{\tilde \kt}
\newcommand{\bcp}{\Bbb C \Bbb P}

 \newcommand{\khx}{K_{hol}(X)}

\newcommand{\khnp}{\kh(\prod_n \bp^1)}
\newcommand{\ktnp}{\kt(\prod_n \bp^1)}
\newcommand{\spp}{ \bz \times SP^\infty(\bp^\infty)}

\title[Holomorphic $K$ -theory]{Holomorphic $K$ - theory, Algebraic
  co-cycles, and Loop Groups }
\author[R.L. Cohen]{Ralph L. Cohen}
 \address{Dept. of Mathematics \\
Stanford University\\
Stanford, California 94305}
 \email[Cohen]
{ralph@math.stanford.edu}
\thanks{The first author was partially supported by a grant from the NSF and a visiting
fellowship from St. Johns College, Cambridge.} 
\author[P. Lima-Filho]{Paulo Lima-Filho}
\address{Department of Mathematics\\
Texas A\&M University \\
College Station, Texas}
\email[Lima-Filho]{plfilho@math.tamu.edu}
\thanks{The second author was partially supported by a grant from the NSF}
\date{\today}
  
\begin{document}

\begin{abstract}
In this paper we   study the ``holomorphic $K$ -theory" of a projective
variety.  This $K$ - theory is defined in terms of the homotopy type of spaces
of holomorphic maps from the variety to Grassmannians and loop
groups. This theory has been introduced in various places such as
\cite{llm-lon},  \cite{friedwalk}, and a related theory was considered
in \cite{karoubi}.    This theory is built out of studying algebraic
bundles over a variety up to ``algebraic equivalence". In this paper
we will give calculations of this theory for ``flag like varieties" 
which include projective spaces, Grassmannians, flag manifolds, and
more general homogeneous spaces, and also give a complete calculation
for symmetric products of projective spaces.  Using  the algebraic
geometric definition of the Chern character studied by the authors in

 \cite{cohenlima}, we will show that there  
 is  a rational isomorphism of graded rings  
  between holomorphic $K$ - theory and the appropriate   
``morphic  cohomology" groups, defined in \cite{friedlaws} in terms of
 algebraic co-cycles in the variety.  In so doing we describe a
 geometric model for  rational morphic  cohomology groups in terms of
 the homotopy type of the space of algebraic maps from the variety to
 the ``symmetrized loop group" $\om U(n)/\Sigma_n$ where the
 symmetric group $\Sigma_n$ acts on $U(n)$ via conjugation. This is
 equivalent to studying algebraic maps to the 
quotient of the infinite Grassmannians $BU(k)$ by a similar symmetric
group action. We then use the Chern character isomorphism   to prove a
conjecture of Friedlander and Walker stating that if one localizes
holomorphic $K$ - theory by inverting the Bott class, then rationally
this is isomorphic to topological $K$ - theory.  Finally this will
 allows us to produce explicit obstructions to periodicity in
 holomorphic $K$ - theory, and show that these obstructions vanish for
 generalized flag manifolds. 
  
\end{abstract}
\maketitle

\section*{Introduction  }

The study of the topology of holomorphic mapping spaces $Hol(X,Y)$,
where $X$ and $Y$ are complex manifolds has been of interest to
topologists and geometers for many years. 
In particular when $Y$ is a Grassmannian or a loop group, the space of
holomorphic maps yields parameter spaces for certain  moduli spaces of
holomorphic bundles (see \cite{ps}, \cite{atiyah}, \cite{cls}).  In
this paper we  study the $K$ -theoretic properties of such holomorphic
mapping spaces.

More specifically, let $X$ be any projective variety,   let $
Gr_n(\bc^M)$ denote the Grassmannian of $n$ - planes in $\bc^M$ (with
its usual structure as a smooth projective variety), and let $\Omega
U(n)$ denote the loop group of the unitary group $U(n)$, with its
structure as an infinite dimensional smooth algebraic variety
(\cite{ps}). We let 
$$
Hol(X; Gr_n(\bc^M)) \quad \text{and} \quad Hol (X; \Omega U(n))
$$ 
denote the spaces of algebraic maps between these varieties
(topologized as  subspaces of the corresponding spaces of continuous
maps, with the compact open topologies).  We use this notation because
if $X$ is smooth these spaces of algebraic maps are precisely the same
as holomorphic maps between the underlying complex manifolds.  The \sl
holomorphic $K$ -theory \rm space $\kh(X)$ is defined to be the
Quillen - Segal group completion of the union of these mapping spaces,
which we write as 
$$
\kh(X) = Hol (X; \bz \times BU)^+ = Hol (X; \Omega U)^+.
$$ 
This group completion process will be described carefully below.  The
holomorphic $K$ - groups will be defined to be the homotopy groups 
$$
\kh^{-q}(X) = \pi_q(\kh(X)).
$$
A variant of this construction was first incidentally
introduced in \cite{llm-lon}, and subsequently developed in
\cite{lmf-mex} where one obtains various connective spectra associated
to an algebraic variety $X$, using spaces of algebraic cycles. The
case of Grassmannians is the one
treated here.  A theory related to   holomorphic $K$ - theory was also studied
by Karoubi in \cite{karoubi} and the construction we use here coincides
with the definition of  ``semi-topological $K$ - theory" studied by Friedlander
and Walker in \cite{friedwalk}.  Indeed their terminology reflects the
fact that for a smooth projective variety $X$ holomorphic $K$ - theory
sits between algebraic $K$ - theory of the associated scheme and the
topological $K$ - theory of  its underlying topological space.  More
precisely, using Morel and Voevodsky's description algebraic $K$
theory  of a smooth variety $X$ (via their work on $\ba^1$ - homotopy 
theory \cite{morvoevod}), Friedlander and Walker showed that there are
natural transformations 
$$
\begin{CD}
\ka (X) @>\alpha >> \kh(X) @>\beta >> \kt(X)
\end{CD}
$$
so that the map $\beta : \kh (X) \to \kt(X)$ is the map induced by
including the holomorphic mapping space $Hol(X; \bz \times BU)$ in the
topological mapping space $Map(X; \bz \times BU)$, 
and where the composition $\beta \circ \alpha : \ka(X) \to \kt(X)$ is
the usual transformation from algebraic $K$ -theory to topological $K$
-theory induced by forgetting the algebraic stucture of a vector
bundle. 

In this paper we calculate the holomorphic $K$ -theory of a large
class of varieties, including ``flag - like varieties", a class that
includes Grassmannians, flag manifolds and more general  
homogeneous spaces.  We also give a complete calculation of the
holomorphic $K$ -theory of  arbitrary symmetric products of projective
spaces.  Since the \sl algebraic \rm $K$ -theory of 
such symmetric product spaces is not in general known, these
calculations should be of interest in their own right.  We then study
the   Chern character for holomorphic $K$ - theory, using the
algebraic geometric description of the   Chern character constructed
by the authors in \cite{cohenlima}. The target  of the Chern character
transformation is the  ``morphic cohomology"
$L^*H^*(X) \otimes \bq$, defined in terms of algebraic co-cycles in $X$
\cite{friedlaws}.  We then prove the following. 

\med
\begin{theorem}\label{E:chiso}
For any projective variety (or appropriate colimit of project
varieties) $X$, the Chern character is a natural transformation 
$$
ch : \kh^{-q}(X) \otimes \bq \la \bigoplus_{k=0}^\infty L^kH^{2k-q}(X)
\otimes \bq 
$$
which is an isomorphism for every $q \geq 0$.  Furthermore it
preserves a natural multiplicative 
structure, so that it is an isomorphism of graded rings.
\end{theorem}

\med
In the proof of this theorem we develop techniques which will yield
the following interesting descriptions of  morphic cohomology that
don't involve the use of higher Chow varieties. 

Consider the following quotient spaces by appropriate actions of the symmetric
groups:

\begin{align}
\Omega U/\Sigma &= \varinjlim_n \Omega U(n)/\Sigma_n, \notag \\
BU/\Sigma &= \varinjlim_{n,m} Gr_m(\bc^{nm})/\Sigma_m \notag\\
SP^\infty(\bcp^\infty) &= \varinjlim_n (\prod_n \bcp^\infty )/\Sigma_n \notag
\end{align}

\med
\begin{theorem}\label{E:symmet}   If
  $X$ is any projective variety, then the Quillen - Segal  group
  completion of the following spaces of algebraic maps 
$$
Mor(X, \Omega U/\Sigma)^+ \quad Mor(X, BU/\Sigma)^+, \quad \text{and} 
\quad Mor(X, SP^\infty (\bcp^\infty))^+
$$
are all rationally homotopy equivalent.  Moreover  their $k^{th}$ -
rational homotopy groups (which we call $\pi_k$) are  isomorphic to 
the rational   morphic cohomology groups
$$
\pi_k \cong \oplus_{p=1}^\infty L^pH^{2p-k}(X) \otimes \bq.
$$
\end{theorem}

\med
Among other things, the relation between morphic cohomology and the
morphism space into the ``symmetrized" loop group allows, using loop
group machinery, a geometric description of these cohomology groups in
terms of a certain moduli space of algebraic  bundles with symmetric  
group action.

We then use   Theorem \ref{E:chiso} to prove the following result about
``Bott periodic holomorphic $K$ - theory".   $\kh^* (X)$ is  a module over
$\kh^*(point)$ in the usual way, and since $\kh^*(point) =
\kt^*(point)$, 
we have a ``Bott class" $b \in \kh^{-2}(point)$.  The module structure
then defines a transformation 
$$
b_* : \kh^{-q}(X) \to \kh^{-q-2}(X).
$$
If $\kh^*(X)[\frac{1}{b}]$ denotes the localization of $\kh^*(X)$
obtained by inverting this operator, 
we will then prove the following rational version of a conjecture of
Friedlander and Walker \cite{friedwalk}. 

\med
\begin{theorem}\label{E:bottlocal}  The map $\beta :
  \kh^*(X)[\frac{1}{b}]\otimes \bq \to \kt^*(X) \otimes \bq $ is an
  isomorphism. 
\end{theorem}

\med
Finally we describe a necessary conditions for the holomorphic $K$ -
theory of a smooth variety to be 
Bott periodic (i.e $\kh(X) \cong \kh(X)[\frac{1}{b}]$) in terms of the
Hodge filtration of its cohomology. We will show that generalized flag
varieties   satisfy this condition and  their holomorphic $K$ - theory
is Bott periodic.  We also give examples of varieties for which these
conditions fail and hence whose holomorphic $K$ - theory is not Bott
periodic. 

\med
This paper is organized as follows.  In section 1 we give the
definition of holomorphic $K$ - theory in terms of loop 
groups and Grassmannians, and prove that the holomorphic $K$ - theory
space, $\kh(X)$ is an infinite loop space.  In section 2 we give a
proof of a result of Friedlander and Walker that $\kh^0(X)$ is the
Grothendieck group of the monoid of algebraic bundles over $X$ modulo
a notion of ``algebraic equivalence".  We prove this theorem here for
the sake of completeness, and also because our proof allows us to
compute the holomorphic $K$ - theory of flag - like varieties, which
we also do in section 2.  In section 3 we identify the equivariant
homotopy type of the holomorphic $K$ - theory space $\kh(\prod_n
\bp^1)$,  where the group action is induced by the permutation action
of the symmetric group $\Sigma_n$. This will allow us to compute the
holomorphic $K$ - theory of symmetric products of projective spaces,
$\kh(SP^m(\bp^n))$.      In section 4 we recall the Chern character 
defined in \cite{cohenlima}   we prove that  
  is an isomorphism  of rational  graded rings(Theorem
  \ref{E:chiso}). In section 5 we prove Theorem \ref{E:symmet} giving
  alternative descriptions of morphic cohomology.  Finally in section
  6 rational maps in the holomorphic $K$ - theory spaces $\kh(X)$ are
  studied, and they are used, together with the Chern character
  isomorphism, to prove Theorem \ref{E:bottlocal} regarding Bott
  periodic  holomorphic $K$ - theory.  
 
The authors would like to thank many of their colleagues for helpful
conversations regarding this work.  They include G. Carlsson,
D. Dugger, E. Friedlander, M. Karoubi, B. Lawson, E.  Lupercio, 
J. Rognes,  and G. Segal.

 
\section{The Holomorphic $K$-theory space}  

In this section we define the holomorphic $K$-theory space $\kh(X)$
for a 
projective variety $X$ and show that it is an infinite loop space.

\med
For the purposes of this paper we let $\om U(n)$ denote the
group of based \sl algebraic \rm loops in the unitary group $U(n)$. That is,
an element of $\om U(n)$ is a map $\gamma : S^1 \to U(n)$ such that $\gamma (1)
= 1$ and $\gamma $ has finite Fourier series expansion.  Namely, $\gamma $ can
be written in the form
$$
\gamma (z) = \sum_{k=-N}^{k=N} A_k z^k
$$
for some $N$, where the $A_k$'s are $n \times n$ matrices.  It is well known
that the inclusion of the group of algebraic loops into the space of all
smooth (or continuous) loops is a homotopy equivalence of infinite
dimensional complex manifolds \cite{ps}.

Let $X$ be a  projective variety.  It was shown by Valli in \cite{valli}
that the holomorphic mapping space $Hol(X, \om U(n))$ has a  $\cc_2$ -
operad structure in the sense of May \cite{may}. Here $\cc_2$ is the
 little 2-dimensional cube operad.  This in particular implies that
the Quillen - Segal group completion, which we denote with the superscript
$+$ (after Quillen's $+$ - construction), $Hol(X, \om U(n))^+$ has the
structure
of a two - fold loop space. (Recall that up to homotopy, the Quillen -
Segal group completion of a topological monoid $A$ is the loop space of the
classifying space, $\Omega BA$.)   By taking the limit over $n$, we define
the holomorphic $K$-theory space to be  the group completion of the
holomorphic mapping
space.

\bg
\begin{definition}
$$
\khx = Hol (X, \om U)^+.
$$
\rm If $A \subset X$ is a    subvariety, we then define the \sl
relative \rm
holomorphic $K$-theory
$$
K_{hol}(X,A)
$$
to be the homotopy fiber of the natural restriction map, $\khx \to
K_{hol}(A).$
\end{definition}

\med
Before we go on we point out certain basic properties of $\khx$.
\smallskip

\noindent{\bf 1.}\ By the geometry of loop groups studied in \cite{ps}
(more specifically the ``Grassmannian model of a loop group") one
knows that every element of the algebraic loop group $\om U(n)$ lies
in a finite dimensional Grassmannian.
When one takes the limit over $n$, it was observed in \cite{cls} that one
has the holomorphic diffeomorphism $\bz \times BU \cong \om U$, where here
$BU$ is given the complex structure as a limit of Grassmannians, and $\om U$
denotes the limit of the \sl algebraic \rm loop groups $\om U(n)$. Thus we
could have replaced $\om U$ by $\bz \times BU$ in the definition of $\khx$.
That is,
we have an equivalent  definition:

\med
\begin{definition}
$$
\kh(X) = Hol(X, \bz \times BU)^+.
$$
\end{definition}

\med
This definition has the conceptual advantage that
$$
\pi_0 \left(Hol(X, BU(n))\right)= lim_{m\to \infty} \pi_0 \left(Hol
  (X, Gr_n(\bc^m))\right) 
$$
where $Gr_n(\bc^m)$ is the Grassmannian of dimension $n$ linear
subspaces of $\bc^m$.  Moreover this set   corresponds to
equivalence classes of rank $n$ holomorphic bundles over $X$
that are embedded (holomorphically) in an $m$ - dimensional trivial
bundle. 
\smallskip

\noindent{\bf 2.}\ It is necessary to take the group completion in our
definition of $\khx$.  For example, the results of \cite{cls} imply
that 
$$
Hol_*(\bp^1, \om U) \cong \coprod_{k=0}^\infty BU(k)
$$
where $Hol_*$ denotes basepoint preserving holomorphic maps. Thus this
holomorphic mapping space is not
  an infinite loop space without group completing.  In fact after we group
complete we obtain

$$
K_{hol}(\bp^1, *) \cong \bz \times BU
$$
and so we have the ``periodicity"  result
$$
K_{hol}(\bp^1, *) \cong K_{hol}(*).
$$
A more general form of
``holomorphic Bott periodicity"  is contained in D. Rowland's Ph.D
thesis \cite{rowland} where it is shown that
$$
K_{hol}(X \times \bp^1, X ) \cong K_{hol}(X)
$$
for any smooth projective variety $X$.  A more general projective
bundle theorem was proved  in \cite{friedwalk}.

\med
We now observe that the two fold loop space mentioned above for holomorphic
$K$-theory can actually be extended to an infinite loop structure.

\med
\begin{proposition} The space $\kh(X) = Hol(X,\bz \times BU)^+$ is an
infinite loop space.
\end{proposition}

\med
\begin{proof}  Let $\cl_*$ be the linear isometries operad.  That is,
$\cl_m$ is the space of linear (complex) isometric embeddings of
$\oplus_m \bc^\infty$ into $\bc^\infty$.  
These spaces are contractible, and the usual operad action
$$
\cl_m \times_{\Sigma_m} (Gr_n(\bc^\infty))^m \la Gr_{nm}(\bc^\infty)
$$
give  holomorphic maps for each $\alpha \in \cl_m$.  
It is then simple to verify that this endows the
holomorphic mapping space 
$$\amalg_n Hol(X, Gr_n(\bc^\infty))$$ with the structure of a $\cl_*$ -
operad space.
Since this is an $E_\infty$ operad in the sense of May \cite{may},
this implies that the group completion, 
$\kh(X) = (\amalg_nHol(X,Gr_n(\bc^\infty))^+$ 
has the structure of an infinite loop space. 
\end{proof}

\med

As is usual, we define the (negative) \sl holomorphic $K$ - groups \rm to
be the homotopy groups of this infinite loop space:

\med
\begin{definition} For $q \geq 0$,
$$
\kh^{-q}(X) = \pi_q(\kh(X)) = \pi_q(Hol(X, \Omega U)^+).
$$
\end{definition}

\med
\noindent{\bf Remarks.}
 
\noindent{\bf a.}\ Notice that as usual, the holomorphic $K$-theory is
a ring.  Namely, the spectrum (in the sense of stable homotopy theory)
corresponding to the infinite loop space $\kh(X)$, is in fact a ring 
spectrum.  The ring structure is induced by tensor product operation on
Grassmannians,
$$
Gr_n(\bc^m)^{\otimes k} \to Gr_{nk}(\bc^{mk}).
$$
We leave it to the reader to check the details that this structure does
indeed induce a ring structure on the holomorphic $K$-theory. Indeed,
this  parallels the well-known fact 
that whenever $\mathbf E$ is a ring spectrum and $X$ is an
arbitrary space, then $\text{Map}(X,{\mathbf E})$ has a natural
structure of ring spectrum, where $\text{Map}(-,-)$ denotes the space
of continuous maps, with the appropriate compact-open, compactly
generated topology; cf. \cite{ may-ring}.

\noindent{\bf b.}\ A variant of this construction was first incidentally
introduced in \cite{llm-lon}, and subsequently developed in
\cite{lmf-mex} where one obtains various connective spectra associated
to an algebraic variety $X$, using spaces of algebraic cycles. The
case of Grassmannians is, up to $\pi_0$ considerations, the one
treated here. A  theory related to  holomorphic $K$ - theory was also
studied by Karoubi in \cite{karoubi}, and the definition given here
coincides with the notion of   ``semi - topological $K$- theory"
introduced and studied by Friedlander and Walker in \cite {friedwalk}.


\section{The holomorphic $K$ - theory of flag varieties and a general
description of $\kh^0(X)$.}

The main goal of this section is to prove the following theorem  which
yields an effective calculation of $\kh^0(X)$, when $X$ is a flag variety.

\med
\begin{theorem}\label{E:flagkh} Let $X$ be a generalized flag variety.
That is, $X$ is a homogeneous space of the form $X = G/P$ where  $G$ is a
complex algebraic group and $P < G$ is a parabolic subgroup.  Then the
natural map from holomorphic $K$ - theory to topological $K$ - theory, 
$$ 
\beta : \kh^0(X) \la \kt^0(X)
$$
is an isomorphism
\end{theorem}

\med
The proof of this theorem involves a comparison of holomorphic $K$ -
theory with algebraic $K$ - theory. As a 
consequence of this comparison we will recover Friedlander and Walker's 
description of $\kh^0(X)$ for any smooth projective variety
$X$ in terms of ``algebraic equivalence classes" of algebraic bundles
\cite{friedwalk}.  We begin by defining this notion of algebraic
equivalence. 

\med  
\med
\begin{definition} Let $X$ be a projective variety (not necessarily
smooth), and
$E_0
\to X$ and
$E_1 \to X$ algebraic bundles.  We say that $E_0$ and $E_1$ are \sl
algebraically equivalent \rm if there exists a constructible, connected
algebraic curve $T$ and an algebraic bundle $E$ over $X
\times T$, so that the restrictions of $E$ to $X \times \{t_o\}$ and $X \times \{t_1\}$ 
are $E_0$ and $E_1$ respectively, for some $t_0, t_1 \in T$.  Here a \sl constructible \rm curve
means a finite union of irreducible algebraic (not necessarily complete) curves in some
projective space.  
\end{definition}

\med
Notice that two algebraically equivalent bundles are \sl topologically \rm
isomorphic, but \sl not \rm necessarily isomorphic as algebraic bundles.

\med
\begin{theorem}\label{E:kh0}  For any  smooth projective algebraic variety $X$, the group
$\kh^0(X)$ is  isomorphic to the Grothendieck group completion of the monoid of algebraic
equivalence classes of  algebraic bundles over $X$.
\end{theorem}

\med

\med
The description of $Mor(X, BU(n))$ given in \cite{cls} provides our
first step in understanding $\kh^0(X)$.

\med
 As in \cite{cls},   if $X$ is a projective variety then we call an algebraic bundle
$E \to X$ \sl embeddable, \rm if there exists an algebraic embedding of $E$ into a trivial bundle:
$E \hookrightarrow X \times \bc^N$ for some large $N$.  Let $\phi : E \hookrightarrow X \times
\bc^N$ be such an embedding. We identify an embedding $\phi$ with the composition $\phi : E \hookrightarrow X
\times \bc^N \hookrightarrow  X \times \bc^{N + M}$, where $\bc^N$ is included in $\bc^{N+M}$ as the first
$N$ coordinates. We think of such an equivalence class of embeddings as an embedding $E \hookrightarrow
X \times \bc^{\infty}$. 
 We refer to the pair $(E,
\phi)$ as an
\sl embedded
\rm algebraic bundle. 

Let $X$ be any projective variety, and let
  $E$ be  a rank $k$   holomorphic bundle over $X$ that is holomorphically embeddable in a
trivial  bundle, define $Hol_E(X, BU(k))$ to be the space of holomorphic maps $\gamma : X
\to BU(k)$ such that
$\gamma^*(\xi_k) \cong E$, where $\xi_k \to BU(k)$ is the universal holomorphic bundle. This is topologized as a subspace
of the continuous mapping space, which is endowed with the compact - open topology.

Let $Aut(E)$ be the gauge group of holomorphic bundle automorphisms of
$E$.  The following lemma identifies the homotopy type of  $Hol_E(X,
BU(k))$ in terms of $Aut(E)$. 

\med
\begin{lemma}\label{E:classif} $Hol_E(X, BU(k))$ is naturally homotopy equivalent to the classifying space
$$
Hol_E(X, BU(k))  \simeq B(Aut(E)).
$$
\end{lemma}

\med
\begin{proof}  As was described in \cite{cls},   elements in  $Hol_E(X, BU(k))$ are in bijective
correspondence to   isomorphism classes of rank $k$ embedded holomorphic bundles, $(\zeta, \phi)$.   By modifying
the embedding $\phi$ via an isomorphism between $\zeta $ and $E$, we see that  $Hol_E(X, BU(k)$ is homeomorphic
to the space of 
 holomorphic embeddings  $\psi : E \hookrightarrow  X \times \bc^\infty$, modulo the action of the
holomorphic automorphism group, $Aut(E)$.  The space of holomorphic embeddings of $E$  in an infinite
dimensional trivial bundle is easily seen to be contractible \cite{cls}, and the action of $Aut(E)$
is clearly free, with local sections.  Again, see \cite{cls} for details.
The lemma follows. 
\end{proof}

\med
\begin{corollary}\label{E:connect}  The space $Hol_E(X; BU(k))$ is connected.
\end{corollary}
\med

Now as above, we say that two embedded algebraic bundles $(E_0, \phi_0)$ and 
$(E_1, \phi_1)$, are \sl path equivalent  \rm if there is topologically embedded bundle 
$(E, \phi)$, over $X \times I$, which gives a path equivalence between $E_0$ and $E_1$, and
over each slice $X \times \{t\}$ is an   embedded algebraic bundle.  Finally, notice that
the set of (algebraic) isomorphism classes of embedded  algebraic
 bundles forms an abelian monoid.  

\med
\begin{lemma}  For any projective algebraic variety $X$, the group
$\kh^0(X)$ is isomorphic to the Grothendieck group completion of the
monoid of path equivalence classes of  embedded algebraic bundles over
$X$.   
\end{lemma}

\med
\begin{proof}  Recall that 
$$
\kh^0(X) = \pi_ 0  \left( \coprod_n Mor(X, BU(n))\right)^+ .
$$
But the set of path components of the Quillen - Segal group completion of a topological $E_\infty$ space is the
Grothendieck group completion of the discrete monoid of path components of the original  $E_\infty$ - space. 
Now as observed  above the morphism space $Hol(X, BU(n))$ is given by configurations of isomorphism
classes of embedded algebraic bundles,  $(E, \phi)$.   Thus $\pi_0(Hol(X, BU(n))$ is the set of path
equivalence classes of such pairs; i.e the set of path equivalence classes of embedded algebraic
bundles of rank $n$.    
We may therefore conclude that 
 $\kh^0(X)$ is the Grothendieck group completion of   the monoid of path equivalence classes
of embeddable algebraic bundles.
\end{proof}

\med

We now strengthen this result as follows.

\med
\begin{lemma} Two embedded bundles $(E_0, \phi_0)$ and $(E_1, \phi_1)$ are path equivalent
if and only if they are algebraically equivalent.
\end{lemma}

\med
\begin{proof}  Let $f : X \times I \to BU(n)$ be the (continuous) classifying map for the
topological bundle $E$ over $X \times I$, which gives the path equivalence between $E_0$ and
$E_1$, and denote $f_0$ and $f_1$ the restrictions of $f$ to $X \times \{0\}$ and $X \times
\{1\}$.  Since $X \times I$ is compact, the image of $f$ is contained in some Grassmannian
$Gr_n(\bc^m) \subset BU(n)$.  It follows that $f_0$ and $f_1$ lie in the same path component
of $Hol(X, Gr_n(\bc^m))$.  Since $Hol(X, Gr_n(\bc^m))$ is a disjoint union of constructible
subsets of the Chow monoid $\cc_{dim X}(X \times Gr_n(\bc^m))$, then $f_0$ and $f_1$
lie in the same connected component of a constructible subset in some projective space.
Using the fact that  any two points in an irreducible algebraic
variety $Y$ lie in some irreducible algebraic curve $C \subset Y$ (see
\cite[p. 56]{mumf}), \rm one concludes that any two points in a connected
constructible subset of projective space lie in a connected constructible
curve. Let $T$ be a connected constructible curve contained in 
$Hol(X,Gr_n(\bc^m))$ and  containing $f_0$ and $f_1$.  
Under the canonical identification 
$Hol (T, Hol(X, Gr_n(\bc^m))) \cong Hol(X \times T, Gr_n(\bc^m))$, one
identifies  the inclusion $$i : T \hookrightarrow Hol(X, Gr_n(\bc^m))$$
with an algebraic morphism $\bar i : X \times T \to Gr_n(\bc^m).$  This
map  classifies the desired embedded bundle $E$ over $X \times T$.

The converse is clear.
\end{proof}

\med
The above two lemmas imply the following.

\med
\begin{proposition}\label{E:kh0emb}  For any projective algebraic variety $X$, (not necessarily smooth), the group
$\kh^0(X)$ is  isomorphic to the Grothendieck group completion of the monoid of algebraic
equivalence classes of embedded algebraic bundles over $X$.
\end{proposition}

\med
Notice that Theorem \ref{E:kh0}  implies that we can remove the ``embedded" condition in the statement of
this proposition.  We will show how that can be done later in this section.   

\med
Recall from the last section that the forgetful map from the category  of colimits of projective
varieties  to the category of topological spaces, induces a map of morphism spaces,
$$
Hol(X; \bz \times BU) \to Map(X; \bz \times BU)
$$
which induces a natural transformation
$$
\beta: \kh(X) \to \kt(X).
$$

\med
\begin{corollary}\label{E:holtop}    Let $X$ be a colimit of projective varieties.  Then the
induced map
$\beta :
\kh^0(X) \to \kt^0(X)$  is induced by sending the class of an embedded bundle to its underlying
topological isomorphism type:
\begin{align}
\beta:\kh^0(X) &\to \kt^0(X) \notag\\
[E,\phi]  &\to [E] \notag
\end{align}
\end{corollary}

\med
In order to approach Theorem 1 we need to understand the relationship between algebraic $K$ - theory,
$\ka^0(X)$, and holomorphic $K$ - theory, $\kh^0(X)$ for $X$ a smooth projective variety.   For such a variety
$\ka^0(X)$ is the Grothendieck group of the exact category of algebraic bundles over $X$.  Roughly speaking
the relationship between algebraic and holomorphic $K$ -theories for a smooth variety is the passage from
isomorphism classes of holomorphic bundles to algebraic equivalence classes of holomorphic bundles.  This relationship
was made precise in \cite{friedwalk}
using the Morel - Voevodsky description of algebraic $K$ - theory of a
smooth scheme  in terms of  an appropriate morphism space.  In particular, recall that
$$
\ka^0(X) = Mor_{\ch((Sm/\bc)_{Nis}}(X, R\Omega B(\sqcup_{n \geq 0} BGL_n(\bc))).
$$
where $\ch((Sm/\bc)_{Nis})$ is the homotopy category of smooth schemes
over $\bc$, using the Nisnevich topology.  See \cite{morvoevod} for details.
In particular   a morphism of projective varieties, $f : X \to Gr_n(\bc^M)$ induces
an element in the above morphism space and hence a class $\langle f \rangle \in \ka^0(X). $  It also induces
a class $[f] \in \pi_0(Hol(X; \bz \times BU)^+ = \kh^0(X)$.  As seen in \cite{friedwalk} this correspondence
extends to give a
 forgetful map from the morphisms in the homotopy category $\ch((Sm/\bc)_{Nis})$ to homotopy classes
of  morphisms in the category of colimits of projective varieties.  This defines a natural transformation
$$
\alpha : \ka^0(X) \to \kh^0(X)
$$
for $X$ a colimit of smooth projective varieties.

\med
\begin{lemma}\label{E:surj}  For $X$ a smooth projective variety the transformation
$$
\alpha : \ka^0(X) \to \kh^0(X)
$$
is surjective.
\end{lemma}

\med
\begin{proof}  As observed above, the set of path components of the Quillen - Segal group completion of
a topological monoid $Y$
is the Grothendieck - group completion of the discrete monoid of path components:
$$
\pi_0(Y^+) = (\pi_0(Y))^+.
$$
Therefore we have that
$\kh^0(X)$ is the Grothendieck group completion of $\pi_0(Hol(X, \bz \times BU)$.  Thus every element  $\gamma
\in \kh^0(X)$ can be written as
$$
\gamma = [f] - [g]
$$
where $f$ and $g$ are holomorphic maps from $X$ to some Grassmannian.  By the above observations
$\gamma = \alpha (\langle f \rangle - \langle g \rangle ).
$
\end{proof}

\med
This lemma and Proposition \ref{E:kh0emb} allow us to prove the following interesting splitting property of
$\kh^0(X)$ which is not immediate from its definition.

\med
\begin{theorem}\label{E: split}  Let $$0 \to [F, \phi_F] \to [E, \phi_E] \to [G, \phi_G] \to 0$$ be a short
exact sequence of embedded holomorphic bundles over a smooth projective variety $X$.  Then in $\kh^0(X)$ we have
the relation
$$
[E,\phi_E] = [F, \phi_F] + [G, \phi_G].
$$
\end{theorem}

\med
\begin{proof}  This follows from  Lemma \ref{E:surj} and the fact that short exact sequences split in
$\ka^0(X)$. \end{proof}

\med
Lemma \ref{E:surj} will also allow us to prove Theorem \ref{E:flagkh} which we now proceed to do.
 We begin with a definition.

\med
\begin{definition}  We say that a smooth projective variety $X$ is \sl{flag - like} \rm if the following
properties hold on its $K$ - theory:
\begin{enumerate}
\item the usual forgetful map
$$
\psi : \ka^0(X) \to \kt^0(X)
$$
is an isomorphism, and 
\item $\ka^0(X)$ is generated (as an abelian group) by embeddable holomorphic bundles.
\end{enumerate}
\end{definition}

\med
\noindent \bf{Remark}: \rm We call such varieties ``flag - like" because generalized flag varieties
(homogeneous spaces $G/P$ as in the statement of Theorem \ref{E:flagkh}) satisfy these conditions.
We now state a strengthening of Theorem \ref{E:flagkh} which we prove.

\med
\begin{theorem}\label{E:flag}  Suppose $X$ is a flag - like smooth projective variety.  Then the homomorphisms
$$
\alpha : \ka^0(X) \to \kh^0(X)
$$
and
$$
\beta: \kh^0(X) \to \kt^0(X)
$$
are isomorphisms of rings.
\end{theorem}

\med
\begin{proof}  Let $X$ be a flag - like smooth projective variety.  Since every embeddable holomorphic bundle
is represented by a holomorphic map $f : X \to Gr_n(\bc^M)$, for some Grassmannian, then property (2) implies
that $\ka^0(X)$ is generated by classes $\langle f \rangle$, where $f$ is such a holomorphic map.  But then
$\psi (\langle f \rangle ) \in \kt^0(X)$ clearly is the class represented by $f$ in $\pi_0(Map(X, \bz \times BU)
= \kt^0(X)$.  But   this means that the map $\psi : \ka^0(X) \to \kt^0(X)$
is given by the composition
$$
\beta \circ \alpha : \ka^0(X) \to \kh^0(X) \to \kt^0(X).
$$
But since $\psi$ is an isomorphism this means $\alpha$ is injective.  But we already saw in corollary
\ref{E:surj}
that $\alpha$ is surjective.  Thus $\alpha$, and therefore $\beta$, are isomorphisms.  Clearly
from their descriptions, $\alpha$ and $\beta$ preserve tensor products, and hence are ring
isomorphisms.
\end{proof}

\med
We now use this result to prove Theorem \ref{E:kh0}.
Let $X$ be a smooth, projective variety and let $e : X \hookrightarrow
\bcp^n$ be a projective embedding.   We begin by describing  a construction which will allow
\sl any \rm holomorophic bundle $E$ over $X$ to be viewed as representing an element of $\kh^0(X)$ (i.e
$E$ does not necessarily have to be embeddable).

So let $E \to X$ be a holomorphic bundle over $X$.     Recall  that by tensoring $E$ with a line bundle of
sufficiently negative Chern class, it  will become embeddable.  (This
is dual to the statement that tensoring a holomorphic bundle over a
smooth projective variety with a line bundle with sufficiently large
Chern class produces holomorphic bundle that is generated by global
sections.)
   So for sufficiently large $k$, the bundle $E \otimes O(-k)$ is
embeddable.  Here $O(-k)$ is the $k$ - fold tensor product of the canonical line bundle $O(-1)$
over $\bcp^n$,  which, by abuse of notation, we identify with its restriction to $X$.
  Now choose a holomorphic embedding
$$
\phi : E \otimes O(-k) \hookrightarrow X \times \bc^N.
$$
Then the pair $(E \otimes O(-k), \phi)$ determines an element of $ \kh^0(X)$.

Now   from Theorem \ref{E:flag} we know that
$\ka^0(\bcp^n) \cong \kh^0(\bcp^n) \cong \kt^0(\bcp^n)$ as rings.
But since  $O(-k) \otimes O(k) = 1  \in  \kt^0(\bcp^n)$, this means that
if
$$
\iota_k = \otimes_k \iota : O(-k) = \otimes_k O(-1) \hookrightarrow \otimes_k \bc^{n+1}
$$
is the canonical embedding, then the pair $(O(-k), \iota_k) $ represents an invertible class
in $\kh^0(\bcp^n)$.  We denote its inverse by  $O(-k)^{-1} \in \kh^0(\bcp^n)$, and, as before,
we use the same notation to denote  its restriction to $\kh^0(X)$.

Write $A(E) =  (E \otimes O(-k), \phi)  \otimes O(-k)^{-1}_e \in \kh^0(X)$.

\med
\begin{proposition}  The assignment to the holomorphic bundle $E$ the class
$$
 A (E) =  [E \otimes O(-k), \phi]  \otimes O(-k)^{-1} \in \kh^0(X)
$$
is well defined, and only depends on the (holomorphic) isomorphism class of $E$.
  \end{proposition}

\med
\begin{proof}  We first verify that given any holomorphic bundle $E \to X$, that $A (E)$
is a well defined element of $\kh^0(X)$.  That is, we need to show that this class is independent
of the choices made in its definition.  More specifically, we need to show that
$$
[E \otimes O(-k), \phi]  \otimes O(-k)^{-1}  =  [E \otimes O(-q), \psi]  \otimes O(-q)^{-1}
$$
for any appropriate choices of $k$, $q$, $\phi$, and $\psi$.  We do this in two steps.

\med
\bf Case 1:  k = q. \rm  In this case it suffices to show that $(E \otimes O(-k), \phi)  $
and $(E \otimes O(-k), \psi) $  lie in the same path component of the morphism
space $Hol(X, BU(n))$, where $n$ is the rank of $E$.  Now using the
notation of  Lemma \ref{E:classif}
 we see that these two elements  both lie in $Hol_{E \otimes O(-k)}(X, BU(n))$,
which, as proved in Corollary \ref{E:connect} there is path connected.

\med
\bf Case 2: General Case: \rm  Suppose without loss of generality that $q > k$.  Then clearly
the classes  $(E \otimes O(-k), \phi)  \otimes O(-k)^{-1} $ and $(E \otimes O(-k) \otimes
O(-(q-k)), \phi \otimes \iota_{q-k})
 \otimes O(-k)^{-1}  \otimes O(-(q-k))^{-1} $ represent the same element of $\kh^0(X)$. But
this latter class is $(E \otimes O(-q), \phi \otimes \iota_{q-k})  \otimes O(-q)^{-1} $
which we know by case 1 represents the same $K$ - theory class as $(E \otimes O(-q), \psi)
\otimes O(-q)^{-1} $.

\med
Thus $A (E)$ is a well defined class in $\kh^0(X)$. Clearly the above arguments also verify
that $A (E)$ only depends on the isomorphism type of $E$.
\end{proof}

\med
Notice that this argument implies
that  $\kh^0$
encodes
\sl all
\rm holomorphic bundles (not just embeddable ones). We will use this to complete the proof of theorem
\ref{E:kh0}.

\med
\begin{proof}  Let $X$ be a smooth projective  variety and let $\ch_X$ denote the Grothendieck group
of monoid of algebraic equivalence classes of holomorphic bundles over $X$.
We show that the correspondence $A$ described in the above theorem induces an isomorphism
$$
\begin{CD}
A : \ch_X @>\cong >> \kh^0(X).
\end{CD}
$$
We first show that $A$ is well defined.  That is, we need to know if $E_0$ and $E_1$ are
algebraically equivalent, then $A(E_0) = A (E_1)$.  So let $E \to X \times T$ be an algebraic
equivalence.  Since $T$ is a curve in projective space, we can find a projective embedding of
the product, $ e : X \times T
\hookrightarrow \bcp^n. $  Now for sufficiently
large $k$,  $E \otimes O(-k)$ is embeddable, and given an embedding $\phi_E$, the pair
$(E \otimes O(-k), \phi_E)$ defines an algebraic equivalence between the embedded bundles
$(E_0 \otimes O(-k), \phi_0)$ and $(E_1 \otimes O(-k), \phi_1)$, where the    $\phi_i$
  are the  appropriate restrictions of the embedding $\phi$.  Thus
$$
[E_0 \otimes O(-k), \phi_0] = [E_1 \otimes O(-k), \phi_1] \in \kh^0(X).
$$ Thus

$$
[E_0 \otimes O(-k), \phi_0]\otimes O(-k)^{-1}   = [E_1 \otimes O(-k), \phi_1] \otimes
O(-k)^{-1}\in \kh^0(X).
$$
But these classes are  $A(E_0)$ and $A (E_1)$. Thus $A : \ch_X \to \kh^0(X)$ is well defined.

Notice also that $A$ is surjective.  This is because, as was seen in the proof of the last theorem,
if $(E, \phi)$ is and embedded holomorphic bundle, then $A (E) = [E, \phi] \in \kh^0(X)$. The essential
point here being that the choice of the embedding $\phi$  does not affect the holomorphic $K$ - theory
class, since the space of such choices is connected.

Finally notice that $A$ is injective.  This is follows from two the two facts:
\begin{enumerate}
\item The classes $[O(-k)^{-1}]$ are units in the ring structure of
$\kh^0(X)$, and 
\item If bundles of the form $E_0 \otimes O(-k)$ and $E_1 \otimes O(-k)$
are algebraically equivalent then the bundles $E_0$ and $E_1$ are
algebraically equivalent.
\end{enumerate}
\end{proof}


\section{The equivariant homotopy type of $\khnp$ and the holomorphic $K$-theory of symmetric products of projective
spaces}

\med
The goal of this section is to completely identify the holomorphic $K$ - theory of symmetric products of projective
spaces, $\kh(SP^n(\bp^m))$.  Since the algebraic $K$ -theory of these spaces is not in general known, this will give
us new information about algebraic bundles over these symmetric product spaces.  These spaces are particularly important
in this paper since, as we will point out below, symmetric products of projective spaces are representing spaces for morphic
cohomology.

Our approach to this question is to study the equivariant homotopy type of $\khnp$, where
the symmetric group $\Sigma_n$ acts on the holomorphic $K$ - theory space $\khnp = Hol(\prod_n\bp^1; \bz \times
BU)^+$ by permuting the coordinates of $\prod_n\bp^1$. It acts on the topological $K$ -theory space $\kt
(\prod_n\bp^1) = Map(\prod_n \bp^1 ; \bz \times BU)$ in the same way.  The main result of this section is t the
following.

\med
\begin{theorem}\label{E:equiv}  The natural map $\beta : \khnp \to \ktnp$ is a $\Sigma_n$ - equivariant
homotopy equivalence.
\end{theorem}

\med
Before we begin the proof of this theorem we observe the following consequences:

\med
\begin{corollary}Let $G < \Sigma_n$ be a subgroup.  Then the induced map on the $K$ - theories
of the orbit spaces,
$$
\alpha : \kh(\prod_n \bp^1/G) \to \kt(\prod_n \bp^1/G)
$$
is a homotopy equivalence.
\end{corollary}

\med
\begin{proof} By Theorem \ref{E:equiv},  $\alpha : \khnp \to \ktnp$ is a $\Sigma_n$ - equivariant
homotopy equivalence.  Therefore it induces a homotopy equivalence on the fixed point sets,
$$\begin{CD} \alpha : \khnp^G @>\simeq >> \ktnp^G \end{CD}.$$ But these fixed point sets are
$\kh(\prod_n \bp^1/G)$ and $ \kt(\prod_n \bp^1/G)$ respectively. \end{proof}

\med
\begin{corollary}\label{E:khcpn}$\alpha : \kh(\bcp^n) \to \kt(\bcp^n)$ is a homotopy equivalence.
\end{corollary}
\begin{proof}  Let $G = \Sigma_n$ in the above corollary. $\prod_n(\bp^1)/\Sigma_n = SP^n(\bp^1) \cong \bcp^{n}$.
\end{proof}

\med
The following example will be important because as seen earlier, symmetric products of projective
spaces form representing spaces for morphic cohomology.

\med
\begin{corollary}\label{E:spcpn} Let $r$ and $k$ be any positive integers.  Then
$$\alpha : \kh(SP^r(\bcp^k)) \to \kt(SP^r(\bcp^k))$$ is a homotopy
equivalence.
\end{corollary}
\begin{proof} Let $G$ be the wreath product $G = \Sigma_r \int
\Sigma_k$ viewed as a subgroup of the symmetric group
$\Sigma_{rk}$. The obtain an identification of orbit spaces
$$
\left(\prod_{rk}\bp^1\right)/\left(\Sigma_r\int \Sigma_k\right) =
SP^r(SP^k(\bp^1)) \cong SP^r(\bcp^k).
$$
Finally, apply the above corollary when $n = rk$.
\end{proof}

Observe that this corollary gives a complete
calculation of the holomorphic $K$- theory of symmetric products of
projective spaces, since their topological $K$ - theory is known.

\med
In order to begin the proof of Theorem \ref{E:equiv} we need to expand
our notion of holomorphic $K$ - theory to include unions of varieties.
So let $A$ and $B$ be subvarieties of $\bcp^n$, then define $Hol(A\cup
B, \bz \times BU)$ to be the space of those continuous maps on $A \cup
B$ that are holomorphic when restricted to $A$ and $B$. This space
still has the action of the little isometry  operad and so we can take
a group completion and define $\kh(A \cup B) = Hol(A\cup B, \bz \times
BU) ^+$.  If $A \cup B$ is connected, then we can define the reduced
holomorphic $K$ - theory as before, $\tilde \kh (A \cup B)    = $ the
homotopy fiber of the restriction map $\kh(A\cup B) \to \kh(x_0)$,
where $x_0 \in A \cap B$.  With this we can now define the holomorphic
$K$ - theory of a smash product of varieties.

\med
\begin{definition}\label{E:smash}
 Let  $X$ and $Y$ be  connected projective  projective varieties with basepoints $x_0$ and $y_0$ respectively.
 We define $\tkh(X \wedge Y)$ to be the homotopy fiber of the restriction map
$$
\rho : \tkh(X \times Y) \to \tkh(X \vee Y)
$$
where $X\vee Y = \{(x,y_0) \}\cup \{(x_0,y)\}  \subset X \times Y$.
\end{definition}

\med
Recall that in topological $K$ -theory, the Bott periodicity theorem can be viewed as saying the Bott
map $\beta : \tkt(X) \to \tkt(X\wedge S^2)$ is a homotopy equivalence for any space $X$.
In \cite{rowland} Rowland studies the holomorphic analogue of this
result. She studies the  Bott map $\beta : \tkh(X) \to \tkh(X
\wedge\bp^1)$ and, using the index of a family of $\bar \partial$
operators, defines a map $\bar \partial : \tkh(X \wedge \bp^1) \to
\tkh (X)$.  Using a refinement of Atiyah's proof of Bott periodicity
\cite{atiyah}, she proves the following.

\med
\begin{theorem}\label{E:rowland} Given any smooth projective variety $X$, the Bott map
$$
\beta : \tkh(X) \to \tkh(X \wedge \bp^1)
$$
is a homotopy equivalence of infinite loop spaces.  Moreover its homotopy inverse is given by the map
$$
\bar \partial : \tkh(X \wedge \bp^1) \to \tkh(X).
$$
\end{theorem}

\med The fact that the Bott map $\beta$ is an isomorphism also follows
from the ``projective bundle theorem"
of Friedlander and Walker \cite{friedwalk} which was proven independently, using different techniques.
This result in the case when $X = S^0$ was proved in \cite{cls}.  The statement in this case is
$$
 \tkh(\bp^1) \simeq \tkh(S^0) = \bz \times BU = \tkt(S^0) \simeq \tkt(S^2).
$$
Combining this with Theorem \ref{E:rowland} (iterated several times) we get the following:

\med
\begin{corollary}\label{E:ksphere} For a positive integer $k$, let $\bigwedge_k \bp^1 = (\bp^1)^{(k)}$ be the $k$ -fold smash
product of $\bp^1$.  Then we have homotopy equivalences
$$
 \tkh( (\bp^1)^{(k)})   \simeq \bz \times BU \simeq \tkt(S^{2k}),
 $$
\end{corollary}

\med
We will use this result to   prove Theorem \ref{E:equiv}.  We actually will prove a splitting result
for $\khnp$ which we now state.

Let $\cs_n$ denote the category whose objects are (unordered) subsets of $\{1, \cdots , n\}$.
Morphisms are inclusions.  Notice that the cardinality of the set of
objects,
$$
|Ob(\cs_n)| = 2^n.
$$
Notice also that the set of objects $Ob(\cs_n)$ has an action of the
symmetric group $\Sigma_n$ induced by the permutation action of
$\Sigma_n $ on $\{1, \cdots , n\}$.

Let $X$ be a space with a basepoint $x_0 \in X$.
For  $\theta  \in  Ob (\cs_n)$,  define
$$
\prod_\theta X = X^\theta  \subset  X^n
$$
by $  X^\theta= \{( x_1, \cdots , x_n)  \quad \text{such that if $j$
is not an element of $\theta$, then $x_j = x_0\in X$}$\}.
Notice that if $\theta$ is a subset of $\{1, \cdots n\}$ of
cardinality $k$, then $X^\theta \cong X^k$.
The smash product $\bigwedge_\theta X = X^{(\theta)}$ is defined
similarly. The following is the splitting theorem that
will allow us to prove Theorem \ref{E:equiv}.

\med
\begin{theorem}\label{E:holsplit} Let $X$ be a smooth projective variety (or a union of smooth projective
varieties).  Then there is a natural $\Sigma_n$ - equivariant homotopy
equivalence 
$$
J : \tkh(  X^n)  \longrightarrow   \prod_{\theta \in Ob(\cs_n)} \tkh(
X^{(\theta)}). 
 $$ where the action of $\Sigma_n$ on the right hand side is induced
by the permutation action of $\Sigma_n$ on the objects $Ob (\cs_n)$. 
\end{theorem}

\begin{proof}  In order to prove this theorem we begin by recalling
the equivariant stable splitting theorem of a product proved in
\cite{CMT}.  An alternate proof of this can be found in
\cite{cohensplit}.  

\med
Given a space $X$ with a basepoint $x_0 \in X$, let $\Sigma^\infty
(X)$ denote the suspension spectrum of $X$.  We refer the reader to
\cite{lewmay} for a discussion of the appropriate category of
equivariant spectra.

\med
\begin{theorem} \label{E:split}There is a natural $\Sigma_n$  equivariant homotopy equivalence of suspension spectra
$$
\begin{CD}
 J :\Sigma^\infty (X^n) @>\simeq >> \Sigma^\infty (
 \bigvee_{\theta
\in Ob(\cs_n)}(X^{(\theta)}))\end{CD}.$$
\end{theorem}

\med
As a corollary of this splitting theorem we get the following splitting of topological $K$ - theory spaces.

\med
\med
\begin{corollary}\label{E:topsplit} There is a $\Sigma_n$ -equivariant homotopy equivalence of topological
$K$ - theory spaces,
$$
J^* : \prod_{\theta \in Ob\cs_n}\tkt(X^{(\theta)}) \to \tkt (X^n).
$$
\end{corollary}

\begin{proof}  Given to spectra $E$ and $F$, let $sMap(E, F)$ be the spectrum consisting of spectrum
maps from $E$ to $F$.  We again refer the reader to \cite{lewmay} for
a discussion of the appropriate
category of spectra.  If $\Omega^\infty $ is the zero space functor from  spectra  to infinite loop spaces,
then $\Omega^\infty(sMap(E, F)) = Map_\infty(\Omega^\infty (E), \Omega^\infty (F)),$
 where $Map_\infty$ refers to the space of infinite loop maps.

Let $bu$ denote the connective topological $K$ - theory spectrum, whose zero space is $\bz \times BU$.
Now Theorem \ref{E:split} yields a $\Sigma_n$ equivariant homotopy equivalence of the mapping spectra,
$$\begin{CD}
J^* : sMap(\Sigma^\infty (\bigvee_{\theta\in Ob(\cs_n)}(X^{(\theta)}), bu) @>\simeq >>
sMap(\Sigma^\infty (X^n), bu),\end{CD}$$and therefore of infinite loop mapping spaces,
$$
\begin{CD}
J^* : Map_\infty(\Omega^\infty\Sigma^\infty (\bigvee_{\theta\in Ob(\cs_n)}(X^{(\theta)}), \bz \times BU) @>\simeq >>
Map_\infty(\Omega^\infty\Sigma^\infty (X^n),  \bz \times BU). \end{CD}$$  But since $\Omega^\infty \Sigma^\infty(Y)$
is, in an appropriate sense, the \sl free \rm infinite loop space generated by a space $Y$, then given any other infinite loop
space $W$, the space of infinite loop maps, $Map_\infty(\Omega^\infty \Sigma^\infty(Y), W)$ is equal to the space
of (ordinary) maps $Map(Y, W)$.  Thus we have a $\Sigma_n$ - equivariant homotopy equivalence of mapping spaces,
$$\begin{CD}
J^* : Map(  (\bigvee_{\theta\in Ob(\cs_n)}(X^{(\theta)}), \bz \times BU) @>\simeq >>
Map( X^n,  \bz \times BU). \end{CD}$$ \end{proof}

\med
\med
Notice that Theorem \ref{E:holsplit} is the holomorphic version of
Corollary \ref{E:topsplit} .    In order to prove this result, we
need to develop a holomorphic version of the arguments used in proving
Theorem \ref {E:topsplit}.  For this we consider the notion of
``holomorphic stable homotopy equivalence", as follows.  Suppose that
$X$ is a smooth projective variety (or union of  varieties) and $E$ is
a spectrum whose zero space is a smooth projective variety (or a union
of such), define
$sHol(\Sigma^\infty(X), E)$ to be the subspace of
$Map_\infty(\Omega^\infty \Sigma^\infty (X), \Omega^\infty (E))  $ consisting of those infinite loop maps $\phi :
\Omega^\infty \Sigma^\infty (X) \to \Omega^\infty (E)$ so that the composition
$$
\begin{CD}
X \hookrightarrow \Omega^\infty \Sigma^\infty (X) @>\phi >> \Omega^\infty (E)
\end{CD}
$$
is holomorphic.  Notice,  for example, that $sHol(\Sigma^\infty(X); bu) = Hol(X, \bz \times BU)$.

Now suppose $X$ and $Y$ are both smooth projective varieties, (or unions of such).

\begin{definition} \label{E:holshe} A map of
suspension spectra,
$\psi : \Sigma^\infty (X) \to \Sigma^\infty(Y)$ is called a \sl holomorphic stable homotopy equivalence, \rm if
the following two conditions are satisfied.
\begin{enumerate}
\item $\psi$ is a homotopy equivalence of spectra.
\notag 
\item If $E$ is any spectrum whose zero space is a smooth projective variety (or a union of such), then the induced
map on mapping spectra, $\psi^* : sMap(\Sigma^\infty(Y), E) \to sMap(\Sigma^\infty (X), E)$ restricts to a map
$$
\psi^*s : Hol(\Sigma^\infty(Y), E) \to sHol(\Sigma^\infty (X), E)
$$
which is a homotopy equivalence. \notag
\end{enumerate}
\end{definition}

\med
With this notion we can complete the proof of Theorem \ref{E:holsplit}.  This requires a proof of Theorem \ref{E:split}
that will respect holomorphic stable homotopy equivalences.  The
version of this theorem given in \cite{cohensplit} will do this.  We
now recall that proof and refer to \cite{cohensplit} for details. 

Let $X$ be a connected space with basepoint $x_0 \in X$.  Let $X_+$ denote $X$ with a disjoint
basepoint, and let $X\vee S^0$ denote the wedge of $X$ with the two point space $S^0$.  Topologically
$X_+$ and $X\vee S^0$ are the same spaces, but their basepoints are in different connected components.
However their suspension spectra $\Sigma^\infty(X_+)$ and $\Sigma^\infty (X \vee S^0)$ are stably homotopy
equivalent spectra with units (i.e via a stable homotopy equivalence $j : \Sigma^\infty (X_+) \simeq \Sigma^\infty(X\vee
S^0)$ that respects the obvious unit maps
$\Sigma^\infty(S^0) \to \Sigma^\infty(X_+)$ and $\Sigma^\infty(S^0) \to \Sigma^\infty(X \vee S^0)$.)
Moreover it is clear that if $X$ is a smooth projective variety then $\Sigma^\infty(X_+)$ and $\Sigma
^\infty(X\vee S^0)$ are holomorphically stably homotopy equivalent in the above sense. Now by taking smash products
$n$ -times of this equivalence, we get a $\Sigma_n$ - equivariant   holomorphic stable homotopy equivalence,
$$
\begin{CD}
J_n : \Sigma^\infty((X_+)^{(n)}) = (\Sigma^\infty((X_+))^{(n)}  @>j^{(n)}>> (\Sigma^\infty(X\vee S^0))^{(n)} =
\Sigma^\infty((X\vee S^0)^{(n)}).
\end{CD}
$$

Now notice that the $n$ - fold smash product
$(X_+)^{(n)}$ is naturally (and $\Sigma_n$ equivariantly) homeomorphic to the cartesian product $(X^n)_+$.   Notice
also  that the
$n$ fold  iterated smash product of
$X \vee S^0$   is $\Sigma_n$ - equivariantly homeomorphic to the wedge of the smash products,
$$
(X\vee S^0)^{(n)} =  \left(\bigvee_{\theta \in Ob (\cs_n)}X^{(\theta)} \right)\vee S^0.
$$
Thus  $J_n$ gives a $\Sigma_n$ - equivariant    stable homotopy equivalence,
$$
\begin{CD}
J_n : \Sigma^\infty(((X^n)_+)  @>\simeq >> \Sigma^\infty( \left(\bigvee_{\theta \in Ob (\cs_n)}X^{(\theta)} \right)\vee S^0)).
 \end{CD}
$$ which gives a proof of Theorem \ref{E:split}.  Moreover when $X$ is a smooth projective variety (or a union of such)
this equivariant stable homotopy equivalence is a holomorphic one.  In particular, given any such $X$, this implies there is a
$\Sigma_n$ equivariant homotopy equivalence
$$
\begin{CD}
 J_n^* : sHol_*(\Sigma^\infty((X^n)_+) ; bu) @>\simeq >> sHol_*(\Sigma^\infty( \left(\bigvee_{\theta \in Ob
(\cs_n)}X^{(\theta)} \right)\vee S^0)) ; bu).
\end{CD}
$$where $sHol_*$ refers to those maps of spectra that preserve the units.  If we remove the units from each of these
mapping spectra we conclude that   we have a $\Sigma_n$ equivariant   homotopy equivalence
$$
 \begin{CD}
 J_n^* : sHol(\Sigma^\infty(X^n) ; bu) @>\simeq >> sHol(\Sigma^\infty \left(\bigvee_{\theta \in Ob
(\cs_n)}X^{(\theta)} \right) ; bu).
\end{CD}
$$
But these spaces are precisely $Hol(X^n, \bz \times BU)$ and $Hol(\bigvee_{\theta \in Ob
(\cs_n)}X^{(\theta)} ;  \bz \times BU) = \prod_{\theta \in Ob
(\cs_n)}Hol (X^{(\theta)} ;  \bz \times BU)$ respectively.  Theorem \ref{E:holsplit} now follows. \end{proof}

\med
We are now in a position to prove   Theorem \ref{E:equiv}.
\begin{proof}  By theorems \ref{E:holsplit} and \ref{E:topsplit} we have the following homotopy commutative
diagram:
$$
\begin{CD}
\tkh((\bp^1)^n) @>J_n^* >\simeq >  \prod_{\theta \in Ob (\cs_n)}\tkh((\bp^1)^{(\theta)})\\
@V\beta VV   @VV\beta V \\
\tkt((\bp^1)^n) @>J_n^* >\simeq >  \prod_{\theta \in Ob (\cs_n)}\tkt((\bp^1)^{(\theta)}).
\end{CD}
$$
Notice that all  the maps in this diagram are $\Sigma_n$ equivariant,  and by the results of theorems \ref{E:holsplit}
and \ref{E:topsplit}  the horizontal maps are $\Sigma_n$ -equivariant homotopy equivalences.  Furthermore,
by  Corollary \ref{E:ksphere} the maps $\tkh((\bp^1)^{(\theta)}) \to \tkt((\bp^1)^{(\theta)})$ are homotopy
equivalences.  Now since the $\Sigma_n$ action on  $\prod_{\theta \in Ob (\cs_n)}\tkh((\bp^1)^{(\theta)})$
and on $\prod_{\theta \in Ob (\cs_n)}\tkt((\bp^1)^{(\theta)})$ is given by permuting the factors according
to the action of $\Sigma_n$ on $Ob( \cs_n)$, this implies that the right hand vertical map
in this diagram,
$\beta : \prod_{\theta \in Ob (\cs_n)}\tkh((\bp^1)^{(\theta)}) \to \prod_{\theta \in Ob (\cs_n)}\tkt((\bp^1)^{(\theta)})$
is a $\Sigma_n$ -equivariant homotopy equivalence.   Hence the left hand vertical map
$$
\beta : \tkh ((\bp^1)^n) \to \tkt((\bp^1)^n)
$$
is also a $\Sigma_n$ - equivariant homotopy equivalence.  This is the
statement of Theorem \ref{E:equiv}. \end{proof}


\section{The Chern character for holomorphic $K$ - theory}

\med
In this section we study  the Chern character for holomorphic
$K$ - theory that was defined by the authors in \cite{cohenlima}.   The values of this Chern character are in the
rational Friedlander - Lawson ``morphic cohomology groups", $L^*H^*(X) \otimes
\bq$.   Our goal is to show that the Chern character  is
gives  an isomorphism $$ch : \kh^{-q}(X) \otimes \bq \cong  \bigoplus_{k=0}^\infty
L^kH^{2k-q}(X)
\otimes
\bq.$$   

Recall the following basic results about the Chern character proved in
\cite{cohenlima}.    

\med
\begin{theorem}\label{E:char} There is a natural transformation of
  functors from the category of colimits of projective varieties to
  algebras over the rational numbers, 
$$
ch : \kh^{-*}(X) \otimes \bq \la \bigoplus_{k=0}^\infty L^kH^{2k-*}(X) \otimes \bq
$$
that satisfies the following properties.
\begin{enumerate}
\item The Chern character is compatible with the Chern character for topological $K$ - theory.
That is, the following diagram commutes:
$$
\begin{CD}
\kh^{-q}(X) \otimes \bq  @>\beta_* >> \kt^{-q}(X) \otimes \bq \\
@Vch VV  @VV ch V \\
\bigoplus_{k=0}^\infty L^k H^{2k-q}(X)\otimes \bq @>>\phi_*>
\bigoplus_{k=0}^\infty H^{2k-q}(X; \bq)
\end{CD}
$$
where $\phi_*$ is the natural transformation from morphic cohomology
to singular cohomology as defined in \cite{friedlaws} 
\item  Let $ch_k : \kh^{-q}(X) \otimes \bq \to L^kH^{2k-q}(X) \otimes
\bq$ be the projection of $ch$ onto the $k^{th}$ factor. Also let $c_k
: \kh^{-q}(X) \to L^kH^{2k-q}(X)$ be the 
$k^{th}$ Chern class defined in \cite[\S 6]{llm-lon}
(see \cite[\S 4]{lmf-mex} for details).  Then
there is a polynomial relation between natural transformations
$$
c_k = k!\ ch_k + p(ch_1, \cdots , ch_{k-1})
$$
where $p(ch_1, \cdots , ch_{k-1})$ is some  polynomial in the first
$k-1$ Chern characters.  
\end{enumerate}
\end{theorem}

As mentioned above  the goal of this section is to prove the following theorem regarding the
Chern character.   

\med
\begin{theorem}\label{E:cherniso}
For every $q \geq 0$, the Chern character for holomorphic $K$ - theory
$$
ch : \kh^{-q}(X)\otimes \bq  \to \bigoplus_{k \geq 0}L^kH^{2k-q}(X) \otimes \bq
$$
is an isomorphism.
\end{theorem}

\med
\begin{proof}  Recall from \cite{friedlaws} that  the suspension theorem
in morphic cohomology implies that morphic cohomology can be represented
by morphisms into spaces of zero cycles in projective spaces.  Since zero
cycles are given by points in symmetric products this can be interpreted
in the following way.  Let $  SP^\infty(\bp^\infty)$ be the infinite
symmetric product of the infinite projective space.     Given a
projective variety $X$, let $Mor(X, \spp)$ denote the colimit of the  the
algebraic morphism spaces $Mor(X; SP^n(\bp^m))$.

\med

\begin{lemma}\label{E:symm} Let $X$ be a colimit of projective varieties.  Then
$$\pi_q(Mor(X; (\spp)^+) \cong \bigoplus_{k \geq 0}L^kH^{2k-q}(X).  $$
\end{lemma}

\med
Similarly, $\bz \times BU$ represents holomorphic $K$ - theory in the sense that
\begin{equation}\label{E:hol}
\pi_q(Mor(X;\bz \times BU)^+) \cong \kh^{-q}(X).
\end{equation}
Thus to prove Theorem \ref{E:cherniso} we will describe a relationship between the representing
spaces $\bz \times SP^\infty(\bp^\infty)$ and $\bz \times BU$.

Using the identification  in Lemma \ref{E:symm},  let
$$ \iota \in \bigoplus_{k = 1}^\infty L^kH^{2k}(\spp) $$
correspond to the class in $  \pi_0((Mor(\spp; \spp)^+)$ represented by the identity
map $id : \spp \to \spp$.

\med
\begin{lemma}\label{E:tau} There exists a unique class $\tau \in \kh^0(\spp)\otimes \bq$ with
Chern character $ch (\tau) = \iota \in \bigoplus_{k = 0}^\infty L^kH^{2k}(\spp) \otimes \bq $.
\end{lemma}

\med
\begin{proof}  By Corollary \ref{E:spcpn}  in  section 3, we know that
for every $k$ and $n$, $\kh(SP^k(\bp^n)) \to \kt (SP^k(\bp^n))$ is a
homotopy equivalence. It follows that by taking limits we have that
$\kh(\spp) \to \kt(\spp) $ is a homotopy equivalence.  But we also
know from \cite{lima} that the natural map
$$\phi : \bigoplus_{k = 0}^\infty L^kH^{2k}(\spp) \otimes \bq  \to
\bigoplus_{k=0}^\infty H^{2k}(\spp; \bq)$$
is an isomorphism. This is true because for the following reasons.
\begin{enumerate}
\item Since products  $\prod_n(\bp^1)$ have   ``algebraic cell
decompositions" in the sense of \cite{lima},
its morphic cohomology and singular cohomology coincide,
$$
\begin{CD}
\phi : L^kH^p(\prod_n(\bp^1) ) @>\cong >> H^p(\prod_n \bp^1).
\end{CD}
$$ 
\item Since both morphic cohomology and singular cohomology admit
transfer maps (\cite {friedlaws}) there is a natural identification of
$L^kH^p(SP^r(\bp^m)\otimes \bq$ and $H^p(SP^r(\bp^m); \bq)$ with the
$\Sigma_r\int\Sigma_m $ invariants in $L^kH^p(\prod_{rm} \bp^1)
\otimes \bq$ and $H^p(\prod_{rm}\bp^1;  \bq)$ respectively.  Since the
natural transformation $\phi : L^kH^p( \prod_{rm} \bp^1) \to
H^p(\prod_{rm}  \bp^1)$ is equivariant,  then we get an induced
isomorphism on the invariants,
$$
\begin{CD}
\phi: L^kH^p(SP^r(\bp^m)) \otimes \bq @>\cong >> H^p(SP^r(\bp^m); \bq).
\end{CD}
$$
\item By taking limits over $r$ and $m$ we conclude that
$$
\phi : L^kH^p(\spp) \otimes \bq \la H^p(\spp; \bq)
$$
is an isomorphism.
\end{enumerate}

Using this isomorphism and  the compatibility of the Chern character maps in holomorphic and topological
$K$ - theories,  to prove this theorem  it is sufficient to prove that there exists
a unique class
$\tau \in \kt^0((\spp) \otimes \bq$ with (topological ) Chern character
$$
ch (\tau) = \iota \in [\spp; \spp] \otimes \bq \cong \oplus_{k=0}^\infty H^{2k}(\spp; \bq)
$$
where $\iota \in [\spp; \spp]$ is the class represented by the identity map.  But this follows because the Chern character
in topological $K$ - theory, $ch : \tkt^0(X) \otimes \bq \to \oplus_{k=1}^\infty H^{2k}(X; \bq)$ is an isomorphism. \end{proof}

\med
We now show how the element $\tau \in \kh^0(\spp) \otimes \bq$ defined in
the above lemma will yield an inverse to the Chern character
transformation.

\med
\begin{theorem}\label{E:trans}  The element $\tau \in \kt^0(\spp) \otimes \bq$ defines  natural
transformations
$$
\tau_* : \bigoplus_{k \geq 0} L^k H^{2k-q}(X) \otimes \bq \la \kh^{-q}(X) \otimes \bq
$$
such that the composition
$$ch \circ \tau_* : \bigoplus_{k \geq 0} L^k H^{2k-q}(X) \otimes \bq \to \kh^{-q}(X) \otimes \bq   \to
\bigoplus_{k \geq 0} L^k H^{2k-q}(X)\otimes \bq$$ is equal to the identity.
\end{theorem}

\med
\begin{proof}    The set of path components of the Quillen - Segal group completion
 of a topological monoid is the Grothendieck group completion of the
discrete monoid of path components.  If we use the notation
$\cm^{\hat{}}$ to mean the Grothendieck group of a discrete monoid
$\cm$, this says that
{\small{$$
\kh^0(\spp) = \pi_0
(Hol(\spp; \bz \times BU)^+) \cong \left(\pi_0(Hol(\spp; \bz \times BU))\right)^{\hat{}} ,
$$ }}
and hence
$$
\kh^0(\spp) \otimes \bq \cong \left(\pi_0(Hol(\spp; \bz \times BU)_{\bq})\right)^{\hat{}} ,
$$ where the subscript $\bq$ denotes the holomorphic mapping space localized at the rationals.
This means that $\tau$ can be represented as a difference of classes,
$$
\tau = [\tau_1] - [\tau_2]
$$
where $\tau_i \in Hol(\spp; \bz \times BU)_\bq$.

Now consider the composition pairing
{\small{
$$
Hol(X; \spp) \times Hol(\spp; \bz \times BU) \to Hol(X; \bz \times BU) \to Hol(X; \bz \times
BU)^+.
$$ }}
which localizes to a pairing
{\small{$$
Hol(X; \spp)_{\bq} \times Hol(\spp; \bz \times BU) _{\bq}\to Hol(X; \bz \times BU)_{\bq} \to Hol(X; \bz
\times BU)^+_{\bq}.
$$}}
Using this pairing, $\tau_1$ and $\tau_2$ each define transformations
$$
\tau_i : Hol(X; \spp)_{\bq} \to Hol(X; \bz \times BU)_{\bq}^+.
$$
Using the fact that $Hol(X; \bz \times BU)_{\bq}^+$ is an infinite loop space, then the subtraction
map is well defined up to homotopy,
$$
\tau_1 - \tau_2 :   Hol(X; \spp)_{\bq} \to Hol(X; \bz \times BU)_{\bq}^+.
$$
We need the following intermediate result about this construction.

\med
\begin{lemma}\label{E:hspace} For any projective variety (or colimit of varieties) $X$, the map
$$
\tau_1 - \tau_2 :   Hol(X; \spp)_{\bq} \to Hol(X; \bz \times BU)_{\bq}^+.
$$ is a  map of $H$ - spaces.
\end{lemma}

\med
\begin{proof} Since the construction of these maps was done at the representing
space level, it is sufficient to verify the claim in the case when $X$ is a point.
That is, we need to verify that the compositions
\begin{equation}\label{E:first}
\begin{CD}
\spp _{\bq}\times \spp_{\bq} @>(\tau_1 - \tau_2) \times (\tau_1 - \tau_2) >>
(\bz \times BU)_{\bq} \times (\bz \times BU)_{\bq} \\
& @>\mu >> (\bz \times BU)_{\bq}
\end{CD}
\end{equation}
and
\begin{equation}\label{E:second}
\begin{CD}
\spp_{\bq} \times \spp_{\bq}  @>\nu >> \spp_{\bq} @>  (\tau_1 - \tau_2)
>> (\bz \times BU)_{\bq}
\end{CD}
\end{equation}
represent the same elements of $\kh^0( \spp \times \spp) \otimes \bq$, where $\mu$ and $\nu$ are
the monoid multiplications in $\bz \times BU$ and $\spp$ respectively.  But by
Corollary \ref{E:spcpn} of the last section, this is the same as $\kt^0((\spp)
\times (\spp)) \otimes \bq$. Now in the topological category, we know
that the class $\tau \in \kh^0(\spp) \otimes \bq$ is the inverse to
the Chern character and hence induces a rational equivalence of $H$ -
spaces
$$
\begin{CD}
\tau : (\spp)_{\bq} @>\simeq >> (\bz \times BU)_{\bq}.
\end{CD}
$$
This implies that  the compositions \ref{E:first} and
\ref{E:second} repesent the same elements of
$\kt^0(\spp) \otimes \bq$, and hence the same elements in
$\kh^0(\spp)\otimes \bq$. \end{proof}
 
\med
Thus the map
$$
\tau_1 - \tau_2 : Hol(X, \spp)_{\bq} \to Hol(X, \bz \times BU)^+_{\bq}
$$
is an $H$ - map  from a $\cc_\infty$ operad spaces (as described in \S 1), to an infinite loop space.   But
any such rational $H$ - map
 extends in a unique manner up to homotopy, to a map of  $H$ - spaces
of their group completions
$$
\tau_1 - \tau_2 : Hol(X. \spp)^+_{\bq} \to Hol(X, \bz \times BU)^+_{\bq}
.$$
This map is natural in the category of colimits of projective varieties $X$. Since any $H$ - map between
rational infinite loop spaces is homotopic to an infinite loop map, this then defines   a natural
transformation of rational infinite loop spaces,

\begin{equation}\label{E:taudef}
\tau  = \tau_1 - \tau_2 : Hol(X. \spp)^+_{\bq} \to Hol(X, \bz \times BU)^+_{\bq}.
\end{equation}

So when we apply homotopy groups  $\tau$ defines  natural transformations

\begin{equation}\label{E:tauhdef}
\tau_* : \bigoplus_{k\geq 0}^\infty L^kH^{2k-q}(X) \otimes \bq  \la \kh^{-q}(X) \otimes \bq.
\end{equation}

Now notice that if we let $X = \spp$ in  (\ref{E:taudef}),  and
$ \iota \in Hol(\spp, \spp)^+_\bq$ be the   class represented by the identity map, then  by
definition, one has that
$$
\tau (\iota) \in Hol(\spp, \bz \times BU)^+_{\bq}
$$
represents the class $[\tau] \in \kh^0(\spp)$ described in Lemma \ref{E:tau}.  Moreover this lemma tells us
that $ch([\tau]) = [\iota] \in \bigoplus_{k\geq 0}^\infty L^kH^{2k}(\spp)$.  Now as in section 4, we view
the Chern character as represented by an element $ch \in Hol(\bz \times BU; \spp)^+_{\bq} $ which is a map of
rational infinite loop spaces, then   this  lemma tells us that the
elements
$$
ch \circ \tau(\iota) \in Hol(\spp ; \spp)^+_{\bq} $$
 and
$$
\iota \in Hol (\spp ;\spp)^+_{\bq}
$$
are both maps of rational infinite loop spaces and lie in the same
path component of $Hol(\spp; \spp)^+_{\bq}$.    But this implies the
$ch \circ \tau$ and $\iota$ define the homotopic natural
transformations of rational infinite loop spaces, 
$$
ch \circ \tau \simeq \iota : Hol(X, \spp)^+_{\bq} \to Hol(X; \spp)^+_{\bq}.
$$
When we apply homotopy groups this means that
$$
ch \circ \tau = id : \bigoplus_{k\geq 0}^\infty L^kH^{2k-q}(X) \otimes
\bq   \la \bigoplus_{k\geq 0}^\infty L^kH^{2k-q}(X) \otimes \bq 
$$
which was the claim in the statement of Theorem \ref{E:trans}. \end{proof}

\med
We now can complete the proof of Theorem \ref{E:cherniso}.  That is we
need to prove that 
$$
ch : \kh^{-q}(X) \otimes \bq \la \bigoplus_{k\geq 0}^\infty
L^kH^{2k-q}(X) \otimes \bq 
$$
is an isomorphism.  By Theorem \ref{E:trans} we know that $ch$ is
surjective.  In order to show that it is injective, we prove the
following:   

\med
\begin{lemma} The  composition of natural transformations $$\tau_* \circ ch : \kh^{-q}(X)\otimes \bq \to
\bigoplus_{k\geq 0}^\infty L^kH^{2k-q}(X) \otimes \bq  \to \kh^{-q}(X) \otimes \bq$$ is the identity.
\end{lemma}

\med
\begin{proof} These transformations are induced on the representing level by maps of rational
infinite loop spaces,
$$
ch : (\bz \times BU)_\bq \to (\spp)_\bq
$$
and
$$
\tau : (\spp )_\bq \to (\bz \times BU)_\bq.
$$
The composition
$$
\tau \circ ch : (\bz \times BU)_{\bq} \to (\bz \times BU)_\bq
$$
represents an element of rational holomorphic $K$ - theory,
$$
[\tau \circ ch ] \in \kh^0(\bz \times BU)_{\bq}.
$$

Now the fact that $\tau $ is an inverse of the Chern character in topological $K$ - theory tells us
that
$$
[\tau \circ ch] = j \in \kt^0(\bz \times BU)_{\bq},
$$
where $j \in \kt^0(\bz \times BU) = \pi_0 (Map(\bz \times BU; \bz \times BU))$ is the class
represented by the identity map.
But according to the results in \S 2, we know
$$
\kh^0(\bz \times BU)_{\bq} \cong \kt^0(\bz \times BU)_{\bq}.
$$  So by the compatibility of the Chern characters in holomorphic and
topological $K$ - theories,  we  conclude that
$$
 [\tau \circ ch] = j \in \kh^0(\bz \times BU)_{\bq}.
$$
This implies that $\tau \circ ch : (\bz \times BU)_\bq \to (\bz \times BU)_\bq $
 and the identity map $id : (\bz \times BU)_\bq \to (\bz \times BU)_\bq $ induce the same
natural transformations
$Hol (X; \bz \times BU)_\bq^+ \to Hol (X; \bz \times BU)_\bq^+$.  Applying homotopy groups
implies that
 $$\tau_* \circ ch : \kh^{-q}(X)\otimes \bq \to \bigoplus_{k\geq
0}^\infty L^kH^{2k-q}(X) \otimes \bq  \to \kh^{-q}(X) \otimes \bq$$ is the identity as claimed. \end{proof}

\med
This lemma implies that $ch : \kh^{-q}(X) \otimes \bq \to  \bigoplus_{k\geq 0}^\infty L^kH^{2k-q}(X)
\otimes \bq $ is injective.  As remarked above this was the last remaining fact to be verified in the proof
of Theorem \ref{E:cherniso}.  \end{proof}

\med
We end this section with a proof that the total Chern class  also gives a rational isomorphism
in every dimenstion.   
  Namely,
recall the Chern classes
$$
c_k : \kh^{-q}(X) \to L^kH^{2k -q}(X)
$$
defined originally in \cite{llm-lon}. Taking the direct sum of these maps gives us the
\sl total Chern class \rm map,
$$
c : \kh^{-q}(X) \to \bigoplus_{k = 0}^\infty L^kH^{2k -q}(X).
$$
We will prove the following result, which was conjectured by
Friedlander and Walker  in \cite{friedwalk}.

\med
\begin{theorem}\label{E:cclass} The total Chern class
$$
c : \kh^{-q}(X) \otimes \bq \to \bigoplus_{k = 0}^\infty L^kH^{2k -q}(X) \otimes \bq
$$
is an isomorphism for all $q \geq 0$.
\end{theorem}

\med
We note that in the case $q = 0$, this theorem was proved in
\cite{friedwalk}.  The proof in general  will follow quickly from our
Theorem \ref{E:cherniso} stating that the total Chern character is a
rational isomorphism.  

\begin{proof} We first prove that the total Chern class
$$
c :\kh^{-q}(X) \otimes \bq   \la \bigoplus_{k=0}^\infty L^kH^{2k - q}(X) \otimes \bq
$$ is injective.  So suppose that for some $\alpha \in \kh^{-q}(X) \otimes \bq $, we have
that $c(\alpha) = 0.$  So each Chern class $c_q(\alpha) = 0$ for $q \geq 0$.  Now
recall from section 4 that in the algebra of operations between
$\kh^{-q}(X) \otimes \bq $ and $\bigoplus_{k=0}^\infty L^kH^{2k - q}(X) \otimes \bq$, that the
that the Chern classes and Chern character are related by a formula of
the form 
\begin{equation}\label{E:classchar}
c_k = k!\ ch_k + p(ch_1, \cdots , ch_{k-1})
\end{equation}
where  $p(ch_1, \cdots , ch_{k-1})$ is some  polynomial in the first
$k-1$ Chern classes.  So since each $c_q(\alpha) = 0$ then an inductive argument using  (\ref{E:classchar}) 
implies that each
$ch_q(\alpha) = 0$.  Thus the total Chern character
$ch (\alpha) = 0$.  But since the total Chern character is an isomorphism (Theorem \ref{E:cherniso}), this
implies that $\alpha = 0 \in  \kh^{-q}(X) \otimes \bq$.  This proves that the total Chern class operation
is injective.

We now prove that $c :\kh^{-q}(X) \otimes \bq   \la \bigoplus_{k=0}^\infty L^kH^{2k - q}(X) \otimes
\bq$ is surjective.  To do this we will  prove that for every $k$ and element $\gamma \in
L^kH^{2k-q}(X) \otimes \bq$ there is a class $\alpha_k \in \kh^{-q}(X)$ with $c_k(\alpha_k) =
\gamma$ and $c_q(\alpha_k) = 0$ for $q \neq k$.   We prove this by induction on $k$.  So assume this
statement is true for $k \leq m-1$, and we now prove it for $k = m$.
Let $\gamma_m \in L^mH^{2m -q}(X) \otimes \bq$.  Since the total Chern character is an
isomorphism, there is an element 
$\alpha_m \in \kh^{-q}(X) \otimes \bq$ with $ch_m(\alpha_m) = \gamma_m$, and $ch_q(\alpha_m) =
0$  for $q \neq m$.  But formula  (\ref{E:classchar}) implies that
$c_q( \alpha_m) = 0$ for $q < m$, and $c_m( \alpha_m) =
\frac{1}{m!}\gamma_m$. Thus the total Chern class has value
$c(m! \alpha_m) = \gamma_m$.
This proves  that the total Chern class is surjective, and therefore
that it is an isomorphism. 
\end{proof}


\section{Stability of rational maps and Bott periodic holomorphic $K$
  - theory} 
 In this section we study the space of rational maps in the morphism
 spaces used to define holomorphic $K$ - theory.  We will show that
 the ``stability  property" for  rational maps in the  morphism space
 $Hol (X, \bz \times BU)$ amounts to the question of whether Bott
 perioidicity holds  in $\kh^*(X)$.  We then use the Chern character
 isomorphism proved in the last section to prove a conjecture of
 Friedlander and Walker \cite{friedwalk} that rationally, Bott
 periodic holomorphic $K$ - theory is isomorphic to topological $K$ -
 theory. (Friedlander and Walker actually  conjectured that this
 statement is true integrally.) 
\bg
Given a projective variety $Y$ with basepoint $y_0 \in Y$ , let
$Hol_{y_0}(\bp^1, Y)$ denote the space of holomorphic (algebraic) maps
$f : \bp^1 \to Y$  satsfying the basepoint condition $f(\infty) =
y_0$.  We refer to this space as the space of based rational maps in 
$Y$.  In \cite{cjs} the ``group completion" of this space of rational
maps  $Hol_{y_0}(\bp^1, Y)^+$ was defined.  This notion of group
completion had    the property that if
$Hol_{y_0}(\bp^1, Y)$ has the structure of a topological monoid, then
$Hol_{y_0}(\bp^1, Y)^+$  is the Quillen - Segal
group completion.  In general $Hol_{y_0}(\bp^1, Y)^+$ was defined to
 be a space of 
limits of ``chains" of rational maps,  topologized using Morse
 theoretic considerations. We refer the reader to \cite{cjs} for
 details. We recall also from that paper the following definition. 

\med
\begin{definition}\label{E:stabil} The space of rational maps in a
  projective variety $Y$  is said to  \sl stabilize, \rm if the group
  completion of the space of rational maps is homotopy equivalent 
to the space of continuous maps,
$$
Hol_{y_0}(\bp^1, Y)^+ \simeq  \Omega^2 Y.
$$
\end{definition}

\med
In \cite{cjs} criteria for when the rational maps in a projective variety
(or symplectic manifold) stabilize were discussed and analyzed.  In
this paper we study the implications in holomorphic $K$ -theory of the
stability of rational maps in the 
 varieties $Hol(X, Gr_n(\bc^M))$, where $X$ is a smooth projective variety,
 and $Gr_n(\bc^M)$ is the Grassmannian of $n$ - dimensional subspaces
 of $\bc^M$. 
(The fact that the space of morphisms from one projective variety to
another is in turn algebraic is well known.  See, for example
\cite{friedwalk2, friedwalk} for discussions about the algebraic
structure of  morphisms between varieties.)  We actually 
study rational maps in $Hol(X;  \bz \times BU)$, which is a colimit of
projective varieties. In fact we will study rational maps in the group
completion $ Hol(X;  \bz \times BU)^+$ by which we mean the group
completion of the relative morphism space, 
$$
Hol_*(\bp^1; Hol(X;  \bz \times BU)^+) = Hol(\bp^1 \times X, \infty
\times X; \bz \times BU)^+. 
$$

\med
\begin{theorem} Let $X$ be a smooth projective variety.  Then the
  space of rational maps in the group completed morphism space 
$Hol(X;\bz \times BU)^+$ stabilizes if and only if the holomorphic $K$
- theory space $\kh(X)$ satisfies Bott periodicity:
$$
\kh(X) \simeq \Omega^2 \kh(X).
$$
\end{theorem}

\med
\begin{proof}  The space of rational maps in the morphism space
  $Hol(X;\bz \times BU)^+$ stabilizes if and only if the group
  completion of its space of rational maps is the two fold loop space, 

\begin{equation}\label{E:rational}
Hol_*(\bp^1; Hol(X; \bz \times BU))^+ \simeq \Omega^2(Hol(X; \bz \times BU)^+).
\end{equation}
But by definition, the left hand side is equal to $Hol(\bp^1 \times X,
\infty \times X; \bz \times BU)^+ = \kh(\bp^1 \times X; \infty \times
X).$  But by Rowland's theorem \cite{rowland} or by the more general
``projective bundle theorem"  proved in  \cite{friedwalk} we know that
the Bott map 
$$
\beta : \kh(X) \to \kh(\bp^1 \times X; \infty \times X)
$$
is a homotopy equivalence.  Combining this with property
\ref{E:rational}, we have that the space of rational maps in the
morphism space $Hol(X;\bz \times BU)^+$ stabilizes if and only if the
following composition is a homotopy equivalence 
\begin{equation}\label{E:bott}
 \begin{CD}
B: \kh(X)  @>\beta > \simeq > \kh(\bp^1 \times X; \infty \times X)\\ 
= Hol_*(\bp^1; Hol(X; \bz \times BU)^+)
 @>>> \Omega^2 Hol(X; \bz \times BU)^+ = \Omega^2 \kh(X).
\end{CD}
\end{equation} \end{proof}

\med
By applying homotopy groups, the Bott map (\ref{E:bott}) 
$B : \kh(X) \to \Omega^2\kh(X)$ defines a homomorphism
$$
B_*:  \kh^{-q}(X) \to \kh^{-q-2}(X)
$$
Let $b \in \kh^{-2}(point)$ be the image under $B_*$of the unit $1 \in
\kh^{0}(point).$ Clearly this class lifts the Bott class in
topological $K$ - theory, $b \in \kt^{-2}(point)$. Observe further
that $B_*:  \kh^{-q}(X) \to \kh^{-q-2}(X)$ is given by multiplication
by the Bott class $b \in \kh^{-2}(point)$, using the module structure
of $\kh^*(X)$ over the ring $\kh^*(point)$. 
The homomorphism $B_*:  \kh^{-q}(X) \to \kh^{-q-2}(X)$ was studied in
\cite{friedwalk} and it was conjectured there that if
$\kh^*(X)[\frac{1}{b}]$ denotes the localization of $\kh^*(X)$
obtained by inverting the Bott class, then one obtains topological $K$
- theory.  We now prove the following rational version of  this conjecture.

\med
\begin{theorem}\label{E:periodic}  Let $X$ be a smooth projective
  variety. Then the map from holomorphic $K$ - theory to topological
$K$ - theory
$\beta : \kh(X) \to \kt(X)$ induces an isomorphism
$$
\begin{CD}
\beta_*:  \kh^*(X)[1/b] \otimes \bq @>\cong >> \kt^*(X) \otimes \bq.
\end{CD}
$$
\end{theorem}

\med
\begin{proof} 
Consider the Chern character defined on the $\kh^{-2}(point) \otimes \bq$  
$$
ch : \kh^{-2}(point) \otimes \bq\to \oplus_{k}L^kH^{2k-2}(point) \otimes \bq.
$$
Now the morphic cohomology of a point is equal to the usual cohomology of a point, $L^kH^{2k-2}(point) = H^{2k-2}(point)$,
so this group is non zero if and only if $k = 1$.  So the Chern character gives an isomorphism
$$
\begin{CD}
ch : \kh^{-2}(point) \otimes \bq @>\cong >> L^1H^0(point) \otimes  \bq\cong \bq.
\end{CD}
$$
Let $s \in L^1H^0(point) \otimes \bq $ be the Chern character of the Bott class, $s = ch(b)$.  Since the Chern
character is an isomorphism, $s \in L^1H^0(point) \otimes \bq \cong \bq$ is a generator.  We use this notation
for the following reason.

Recall the operation in morphic cohomology
$S : L^kH^{q}(X) \to L^{k+1}H^q(X)$ defined in \cite{friedlaws}.
Using the fact that $L^*H^*(X)$ is a module over $L^*H^*(point)$ (using the ``join" multiplication
in morphic cohomology), then this operation is given by multiplication by a generator
of $L^1H^0(point) = \bz$.    Therefore up  to a rational multiple, this operation on rational morphic cohomology,
$ S : L^kH^{q}(X) \otimes \bq\to L^{k+1}H^q(X) \otimes \bq$, is given by multiplication by the element
$s = ch(b) \in L^1H^0(point) \otimes \bq$.

In \cite{friedlaws} it was shown that the natural map from morphic
cohomology to singular cohomology
$\phi : L^kH^q(X) \to H^q(X)$ makes the following diagram commute:
\begin{equation}\label{E:sphi}
\begin{CD}
L^kH^q(X) @>S >> L^{k+1}H^q(X) \\
@V\phi VV    @VV\phi V\\
H^q(X) @> = >> H^q(X).
\end{CD}
\end{equation}

It also follows from the ``Poincare duality theorem"  proved in
\cite{friedlaws-dual} that if $X$ is an $n$ - dimensional smooth
variety,  then $L^sH^q(X)  = H^q(X)$ for $s \geq n$. Furthermore  for
$k < n$ $\phi : L^kH^q(X) \to H^q(X)$ factors as the composition
\begin{equation}\label{E:scomp}
\begin{CD}
\phi: L^kH^q(X) @>S >> L^{k+1}H^q(X) @>S >> \cdots @>S >> L^nH^q(X) = H^q(X)
\end{CD}
\end{equation}

Let $L^*H^q(X)[1/S]$ denote the localization of $L^*H^q(X)$ obtained by inverting
 the transformation $S : L^*H^q(X) \to L^{*+1}H^q(X)$.  Specifically
$$
\begin{CD}
L^*H^q(X)[1/S] = \varinjlim \{ L^*H^q(X) @>S >> L^{*+1}H^q(X) @>S >> \cdots \}
\end{CD}
$$

 Then (\ref{E:sphi}) and (\ref{E:scomp})   imply we have an isomorphism with singular cohomology,
\begin{equation} 
\begin{CD}
\phi : L^*H^q(X)[1/S] @>\cong >> H^q(X).
\end{CD}
\end{equation}
Again, since rationally the $S$ operation is, up to multiplication by a nonzero rational number, given by multiplication
by $s \in L^1H^0(point) \otimes \bq$, we can all conclude that when rational  morphic cohomology is localized by inverting
$s$,  we have an isomorphism with singular rational cohomology,
\begin{equation}\label{E:slocal}
\begin{CD}
\phi : L^*H^q(X; \bq)[1/s]  @>\cong >> H^q(X; \bq).
\end{CD}
\end{equation}

Now since the Chern character isomorphism $ch :\kh^{-q}(X) \otimes \bq  \to\oplus_{k=0}^\infty L^kH^{2k-q}(X) \otimes
\bq$ is an isomorphism of rings,  then the following diagram commutes:   
\begin{equation}\label{E:bottS}
\begin{CD}
\kh^{-q}(X) \otimes \bq  @>\cdot b >> \kh^{-q-2}(X)  \otimes \bq\\
@V ch V \cong V   @V\cong V ch V \\
\bigoplus_{k=0}^\infty L^kH^{2k-q}(X) \otimes \bq @>> \cdot s > \bigoplus_{k=0}^\infty L^{k+1}H^{2k-q}(X)
\otimes \bq.
\end{CD}
\end{equation} where the top horizontal map is multiplication by the Bott class $b \in \kh^{-2}(point)$, and the bottom
horizontal map is multiplication by $s = ch(b) \in L^1H^0(point) \otimes \bq$.

Moreover since the Chern character  in holomorphic $K$ - theory and and that  for topological $K$ - theory
are compatible, this means we get a commutative diagram:
$$
\begin{CD}
\kh^{-q}(X)[1/b] \otimes \bq @>\beta >> \kt^{-q}(X) \otimes \bq\\
@V ch VV   @VV ch V \\
\bigoplus_{k=0}^\infty L^*H^{2k-q}(X; \bq)   [1/s] @>> \phi  > \bigoplus_{k=0}^\infty H^{2k-q}(X; \bq).
\end{CD}
$$
By (\ref{E:slocal})     we know that the bottom horizontal map is an isomorphism.  Moreover by
Theorem \ref{E:cclass}  the left hand vertical map is an isomorphism.   
 Of course the right hand vertical map is also a rational isomorphsim.
Hence the top horizontal map is a rational isomorphism,
$$
\begin{CD}
\beta _* : \kh^{-q}(X)[1/b] \otimes \bq @>\cong>> \kt^{-q}(X) \otimes \bq.
\end{CD}
$$ \end{proof}

\med

In most of the  calculations  of $\kh(X)$ done so far we have seen examples of when
$\kh(X) \cong \kt(X)$.  In particular in these examples the holomorphic $K$ - theory
is periodic, $\kh^*(X) \cong \kh^*(X)[\frac{1}{b}]$. As we have seen from Theorem \ref{E:periodic},
these two conditions are rationally equivalent.
We end by using the above results to give a necessary condition for the holomorphic $K$ - theory to be
Bott periodic, and use it to describe examples where periodicity fails, and therefore provide examples
that have distinct holomorphic and topological $K$ - theories.

\med
\begin{theorem}Let $X$ be a smooth projective variety. Then if $\kh^*(X)\otimes \bq \cong
\kh^*(X)[\frac{1}{b}]$ (or equivalently $\kh^*(X)\otimes \bq
\cong \kt^*(X)\otimes \bq$), then in the Hodge filtration of its cohomology we have
$$
H^{k,k}(X; \bc) \cong H^{2k}(X; \bc)
$$
for every $k \geq 0$.
\end{theorem}

\med
\begin{proof} Consider the commutative diagram involving the total Chern character
\begin{equation}
\begin{CD}
\kh^0(X)\otimes \bc @>\beta_* >> \kt^0(X)\otimes \bc \\
@Vch VV    @VV ch V  \\
\bigoplus_{k \geq 0}L^kH^{2k}(X) \otimes \bc  @>\phi >> \bigoplus_{k \geq 0} H^{2k}(X; \bc)
\end{CD}
\end{equation}

By Theorem \ref{E:periodic}, if $\kh^*(X)$ is Bott periodic, then the top horizontal map
$\beta : \kh^0(X) \otimes \bc \to \kt^0(X)\otimes \bc$ is an isomorphism.  But by theorem
\ref{E:cherniso} we know that the two vertical maps in this diagram are isomorphisms.  Thus if
$\kh^*(X)$ is Bott periodic,  then the bottom horizontal map in this diagram is an isomorphism.
That is,    $$\phi : L^kH^{2k}(X) \otimes \bc \la H^{2k}(X; \bc)$$ is an isomorphism, for every $k
\geq 0$. But as is shown in \cite{friedlaws}, $L^kH^{2k}(X)  \cong
\ca_k(X)$, where $\ca_k(X)$ is the space of algebraic $k$ - cycles in
$X$ up to algebraic (or homological) equivalence. Moreover the image
of $\phi : L^kH^{2k}(X) \otimes \bc \to H^{2k}(X; \bc)$ is the image
of the natural map induced by including algebraic cycles in all
cycles, $\ca_k \otimes \bc \to H^{2k}(X; \bc)$, which lies in Hodge
filtration $H^{k,k}(X; \bc) \subset H^{2k}(X; \bc)$. Thus
$\phi : L^kH^{2k}(X) \otimes \bc\to H^{2k}(X; \bc)$
is an isomorphism implies that the composition
$$\ca_k(X) \otimes \bc  \to H^{k,k}(X; \bc) \subset H^{2k}(X; \bc)
$$
is an isomorphism.  In particular this means that $H^{k,k}(X; \bc) = H^{2k}(X; \bc)$.
\end{proof}

\med
We end by noting that for a flag manifold $X$, we know by Theorem \ref{E:flagkh} that $\kh^0(X) \cong
\kt^0(X)$, and  indeed $H^{p,p}(X;\bc) \cong H^{2p}(X; \bc)$.  However in general  this
theorem tells us that if have a variety
$X$ having nonzero
$H^{p,q}(X; \bc)$ for some $p \neq q$, then $\kh^*(X)$ is \sl not \rm Bott periodic,
and in particular is distinct from topological $K$ - theory.  Certainly abelian varieties
of dimension $\geq 2$ are examples of such varieties.

\end{document}